\documentclass[namedreferences]{kluwer}

\usepackage{epsfig}
\usepackage{wasysym}

\newtheorem{theorem}{Theorem}
\newtheorem{prop}{Proposition}
\newtheorem{lemma}{Lemma}
\newtheorem{osservazione}{Remark}
\begin{document}

\begin{article}
\begin{opening}
\title{On the Construction of Some Buchsbaum Varieties and the Hilbert Scheme of
Elliptic  Scrolls in $\dP^{5}$}

\author{DOLORES \surname{BAZAN}\email{bazan@mathsun1.univ.trieste.it}}
\institute{Via degli Orti 10, Gradisca (GO), Italy}
\author{EMILIA \surname{MEZZETTI}\email{mezzette@univ.trieste.it}\thanks{supported by MURST and GNSAGA, member
of EAGER}}
\institute{Universit\`a degli Studi di Trieste, Dipartimento di Scienze Matematiche, Via Valerio 12/1, 
34127 Trieste, Italy}

\runningtitle{BUCHSBAUM VARIETIES AND ELLIPTIC SCROLLS}
\runningauthor{BAZAN -- MEZZETTI}

\begin{abstract}
 We study the degeneracy loci of 
general bundle morphisms of the form 
$\mathcal{O}^{\oplus m}_{\dP^{n}}\longrightarrow \Omega_{\dP^{n}}(2)$, also from 
the point of view of the classical geometrical interpretation of the sections of $\Omega_{\dP
^{n}}(2)$ as linear line complexes in $\dP^{n}$. We consider in particular the case of 
$\dP^{5}$ with $m=2, 3$. For $n=5$ and $m=3$ we give an explicit description 
of the Hilbert scheme $\mathcal{H}$ of elliptic normal scrolls in $\dP^{5}$, by defining a 
natural rational map $\rho :\mathbb G(2,14)\dashrightarrow\mathcal{H}\:$, which results to be 
dominant with general fibre of degree four.
\end{abstract}


\classification{Mathematics Subject Classification (1991)}{14F05, 14H52, 14J10}

\end{opening}

\section{Introduction}

Degeneracy loci of general bundle morphisms of the form 
\[ \mathcal{O}^{\oplus m}_{\dP^{n}}\longrightarrow \Omega_{\dP^{n}}(2) \] have been 
studied by several authors, for example M.C. Chang (\cite{Chang}) and G. Ottaviani (\cite{Ot}). We have performed 
here a study of these loci in the general case, computing in particular a locally free resolution 
of their ideals, the cohomology groups of their ideal sheaf and their degree. They result to be 
all arithmetically Buchsbaum varieties, with interesting geometrical properties. This can be seen 
by the classical geometrical interpretation of the sections of $\Omega_{\dP^{n}}(2)$ as 
linear line complexes in $\dP^{n}$. For example, if $n$ is odd, they are scrolls, while, 
if $n$ is even, they are unirational varieties.

\medskip

 We have considered then the particular case of $\dP^{5}$ with $m=2,3$.
If $m=2$, the degeneracy locus $X$ is the union of three lines. We have analyzed the 
possible configurations of these lines, obtaining the result that $X$ has the expected dimension 
$1$ if and only if $X$ is exactly the union of three skew lines generating $\dP^{5}$, i.e. 
without any common secant line.

\medskip

This has been applied to the case $m=3$. Here, if the degeneracy locus $X$ has the expected 
dimension $2$, $X$ is an elliptic normal scroll of degree six. It is known that the Hilbert scheme 
$\mathcal{H}$ of elliptic normal scrolls in $\dP^{5}$ has dimension $36$ (see \cite{Ionescu}), which 
is equal to the dimension of $\mathbb{G}(2,14)$, the Grassmannian parametrizing maps 
$\mathcal{O}^{3}_{\dP^{5}}\longrightarrow\Omega_{\dP^{5}}(2)$. Moreover there is 
an open subset in $\mathcal{H}$ corresponding to elliptic scrolls with only two distinct 
unisecant cubics. There is a natural rational map $\rho :\mathbb{G}(2,14)\dashrightarrow
\mathcal{H}$: we have studied this map, in particular to understand if it is dominant and what 
are its fibres. This problem had been tackled classically by G. Fano in \cite{Fano}. Our result, which 
was already known to Fano, is Theorem \ref{fano}, saying that if $X$ is an elliptic normal scroll in 
$\dP^{5}$ with only two distinct unisecant cubics, then there are exactly four general 
nets of linear line complexes having $X$ as singular surface. So the map $\rho$ is dominant 
with general fibre of degree four.

\medskip

If $X$ is an elliptic normal scroll with only one unisecant cubic or with a $1$-dimensional 
family of unisecant cubics, this result remains still valid, but \emph{one} of the four nets having 
$X$ as singular surface is not general, i.e. it contains some special complexes of second type. 
In this situation, the singular locus of such a net does not have the expected codimension $3$, but $2$: 
indeed it is  the union of $X$ with the singular $3$-space of the special complexes of second 
type contained in the net.

\medskip
      
 We will denote by $\dP^{n}$ the projective space of dimension $n$ over an algebraically 
closed field $K$ of characteristic $0$. If $\cal F$ is a sheaf on $\dP^n$, the direct sum of $m$
copies of $\cal F$ will be denoted by ${\cal F}^{\oplus m}$ or else by $m\cal F$.


\section{Degeneracy locus of a morphism 
${\varphi}: \mathcal{O}^{\oplus m}_{\dP^{n}}\longrightarrow
{\Omega}_{\dP^{n}}(2)$}


 Let $\mathcal{O}_{\dP^{n}}$ be the sheaf of  regular functions on 
$\dP^{n}$ and $\Omega_{\dP^{n}}$ be the cotangent sheaf. $\Omega_{\dP^{n}}$ 
is locally free of rank $n$. We recall that the cohomology groups $H^{0}(\Omega_{\dP^{n}}(k))$ 
are zero if $k<2$, while $\Omega_{\dP^{n}}(2)$ is generated by its global sections and 
$\dim H^{0}(\Omega_{\dP^{n}}(2))={{n+1}\choose {2}}$.

A morphism $\varphi:\mathcal{O}^{\oplus m}_{\dP^{n}}\longrightarrow\Omega_{\dP^{n}}(2)$
 is assigned by giving $m$ global sections of $\Omega_{\dP^{n}}(2)$. If $m\leqslant n$, 
then $m$ general sections of $\Omega_{\dP^{n}}(2)$ are linearly independent, so a \emph{general} 
morphism $\varphi$ is generically injective.

We want to study such a morphism $\varphi$ for $m\leqslant n$ and, in particular, 
its degeneracy locus $X$. From Bertini type theorems we have that, if non-empty, $X$ has the 
expected codimension $n-m+1$, that is $\dim X = m-1$. Furthermore the singular locus of $X$, 
$Sing \ X$, has dimension at most $2m-n-4$. In particular $X$ is smooth if \( m<\frac{n+4}{2} \).

\medskip

Eagon-Northcott's theorem (see for instance \cite{GP}) gives a locally free  resolution of 
$\mathcal{O}_{X}$. Indeed, $\varphi$ is general, therefore the Eagon-Northcott's complexes 
associated to the dual morphism $\varphi ^{\ast}$ are exact and have length equal to the 
codimension of $X$. The first of such complexes can be written as follows:

\begin{eqnarray}\label{EN}
\ 0 &\rightarrow\wedge^{n}(\mathcal{T}_{\dP^{n}}(-2))\otimes S^{n-m}(\mathcal{O}^
{\oplus m}_{\dP^{n}})&\rightarrow\wedge^{n-1}(\mathcal{T}_{\dP^{n}}(-2))\otimes 
S^{n-m-1}(\mathcal{O}^{\oplus m}_{\dP^{n}})\rightarrow  \nonumber\\  
\cdots & \rightarrow
\wedge^{m}(\mathcal{T}_{\dP^{n}}(-2))&\rightarrow\mathcal{O}_{\dP^{n}}\rightarrow
\mathcal{O}_{X}\rightarrow 0
\end{eqnarray}

\noindent with natural maps induced by $\varphi$.
\medskip 

Furthermore the class of $X$ in the Chow ring of $\dP^{n}$ is the top Chern class 
$c_{n-m+1}({\hbox{coker}}\ \varphi)$; in particular we can compute the degree of $X$, which is: 

\[
\deg X =\sum_{i=0}^{n-m+1}(-1)^{i}{n-i\choose m-1}.
\]                   

\noindent We recall the following isomorphisms: 

\[ S^{k}(\mathcal{O}^{\oplus m}_{\dP^{n}})\simeq{m+k-1\choose k}\mathcal{O}_{\dP^{n}}\]

\[\wedge^{j}\mathcal{T}_{\dP^{n}}\simeq (\wedge^{n-j}\mathcal{T}_{\dP^{n}})\otimes 
\wedge^{n}\mathcal{T}_{\dP^{n}}\simeq \Omega^{n-j}_{\dP^{n}}(n+1)\] 

\noindent Substituting in (\ref{EN}) and twisting by $n-1$, we obtain at last :

\begin{eqnarray}\label{EN2} 
0&\rightarrow {n-1\choose m-1}\mathcal{O}_{\dP^{n}}\rightarrow {n-2\choose m-1}
\Omega_{\dP^{n}}(2)\rightarrow {n-3\choose m-1}\Omega^{2}_{\dP^{n}}(4)\rightarrow
\cdots \\ 
\cdots& \rightarrow{n-j-1\choose m-1}\Omega^{j}_{\dP^{n}}(2j)\rightarrow\cdots
\rightarrow\Omega^{n-m}_{\dP^{n}}(2n-2m)\rightarrow\mathcal{I}_{X}(n-1)\rightarrow 0 \nonumber
\end{eqnarray}

\noindent where $\mathcal{I}_{X}$ is the ideal sheaf of $X$. From this resolution, we can compute 
the dimensions of the cohomology groups of $\mathcal{I}_{X}$.  

\begin{prop} The cohomology groups $H^{i}(\mathcal{I}_{X}(p))$ for $i>0$ are all zero, 
except the following ones:

\noindent$H^{1}(\mathcal{I}_{X}(m-2)), H^{3}(\mathcal{I}_{X}(m-4)), \ldots ,H^{n-m}(\mathcal{I}_{X}(2m-n-1))$ if
$n-m$ is  odd,
$H^{2}(\mathcal{I}_{X}(m-3)), H^{4}(\mathcal{I}_{X}(m-5)),\ldots ,H^{n-m}(\mathcal{I}_{X}(2m-n-1))$ if $n-m$ is
even. 

\noindent In particular, X is arithmetically Buchsbaum.
\end{prop}

\emph{Proof}. The dimensions of the groups $H^{i}(\mathcal{I}_{X}(p))$ can be computed from (\ref{EN2}). 
In particular, the multiplication maps \[ H^{i}(\mathcal{I}_{X}(p))\to H^{i}(\mathcal{I}_{X}(p+1))\] 
are all zero for $i>0$ and the following condition of St\"uckrad-Vogel (\cite{SV}, see also \cite{Chang}), which 
ensures that $X$ is arithmetically Buchsbaum, is fulfilled: if $H^{i}(\mathcal{I}_{X}(p))\neq 0$ and 
$H^{j}(\mathcal{I}_{X}(q))\neq 0$ with $i<j$, then $j-i\neq p-q-1$. \
\ $\Box$ 
\medskip

We recall now a very nice geometrical interpretation of a morphism $\varphi$ as above (see \cite{Ot}).

From the Euler sequence for $\dP^{n}=\dP(V)$ twisted by two:
\[  0\rightarrow\Omega_{\dP^{n}}(2)\rightarrow\mathcal{O}_{\dP^{n}}(1)^{\oplus (n+1)}
\rightarrow\mathcal{O}_{\dP^{n}}(2)\rightarrow 0
\] 

\noindent taking global sections and noting that $H^{0}(\mathcal{O}_{\dP^{n}}(1)^{\oplus(n+1)})
\simeq V^{\ast}\times V^{\ast}$ and $H^{0}(\mathcal{O}_{\dP^{n}}(2))\simeq S^{2}V^{\ast}$, 
we obtain $H^{0}(\Omega_{\dP^{n}}(2))\simeq (\wedge^{2}V)^{\ast}$.

 This allows to interpret a global section of $\Omega_{\dP^{n}}(2)$ as a bilinear 
alternating form on $V$, or as a skew-symmetric matrix of type $(n+1)\times(n+1)$ with entries 
in the base field. So a general morphism $\varphi$ is assigned by giving $m$ general skew-symmetric 
matrices $A_{1},\ldots ,A_{m}$ and the corresponding degeneracy locus $X$ in $\dP^{n}$ 
is the set of the points $P$ such that:

\begin{varequation}{3} 
(\lambda_{1}A_{1}+\cdots+\lambda_{m}A_{m})[P]=0  
\end{varequation} 

\noindent for some $(\lambda_{1},\ldots ,\lambda_{m})\neq(0,\ldots,0)$, where $[P]$ denotes the column matrix 
of the coordinates of $P$.

 Since a skew-symmetric matrix always is of even rank, to study equation (3) we have 
to distinguish two cases:
\begin{enumerate}[00]
\item if $n$ is even, for every $m$-tuple $\lambda_{1},\ldots ,\lambda_{m}$ equation (3) 
has at least one solution $P\in X$; in this case $X$ is an unirational variety pa\-ra\-me\-tri\-zed by 
$\dP^{m-1}$;

\item if $n$ is odd, then the vanishing of the Pfaffian of the matrix 
$\lambda_{1}A_{1}+\cdots +\lambda_{m}A_{m}$ defines a hypersurface $Z$ of degree \( \frac{n+1}{2} \) 
in $\dP^{m-1}$, where $\lambda_{1},\ldots ,\lambda_{m}$ are homogeneus coordinates. Furthermore, 
if $\varphi$ is \emph{general}, for a fixed point $[\lambda]\in Z$, we find a line of solutions 
of (3) in $X$, so that $X$ is a scroll over $Z$.
\end{enumerate}

\par \noindent In particular, from these observations we can conclude that $X$ is always non-empty.

To see another interpretation of global sections of $\Omega_{\dP^{n}}(2)$, 
let us consider $\mathbb{G}(1,n)\subseteq \dP(\wedge^{2}V)$, the Grassmannian of lines 
in $\dP^{n}$, embedded in $\dP(\wedge^{2}V)$ via the Pl\"ucker map. The dual space 
$\dP(\wedge^{2}V^{\ast})$ parametrizes hyperplane sections of $\mathbb{G}(1,n)$ or, in the 
old terminology, the \emph{linear line complexes} in $\dP^{n}$. Such a linear complex $\Gamma$ 
is represented by a linear equation in the Pl\"ucker coordinates $p_{ij}$: $\;\sum_{0\leqslant i<j 
\leqslant n}a_{ij}p_{ij}=0$. We can associate to it a skew-symmetric matrix $A=(a_{ij})$ of order $n+1$. 
A point $P\in \dP^{n}$ is called a \emph{centre} of $\Gamma$ if all lines through $P$ belong 
to $\Gamma$. The space $\dP(\ker(A))$ results to be the set of centres of $\Gamma$: it is 
called the \emph{singular space} of $\Gamma$.

Once again we have to distinguish the following two cases:

\begin{enumerate}[00]
\item if $n$ is even, then every linear complex $\Gamma$ possesses at least one centre. 
$\Gamma$ is said to be \emph{special} if its singular space is at least a plane. A special $\Gamma$ 
corresponds to a hyperplane section of $\mathbb{G}(1,n)$ with a tangent hyperplane or, equivalently, 
to a point of $\check{\mathbb{G}}(1,n)$, the dual variety of $\mathbb{G}(1,n)$. It is known that, 
$n$ being even, the dual Grassmannian $\check{\mathbb{G}}(1,n)$ has codimension three in 
$\dP(\wedge ^{2}V^{\ast})$.

\item If $n$ is odd, then a general linear complex $\Gamma$ does not have any centre, whereas it 
is said to be \emph{special} if its singular space is at least a line. As above, if $\Gamma$ is special, 
it corresponds to a tangent hyperplane section of $\mathbb{G}(1,n)$ or, which is the same, to a point 
of $\check{\mathbb{G}}(1,n)$ which in this case is a hypersurface in $\dP(\wedge ^{2}V^{\ast})$ 
of degree \( \frac{n+1}{2} \).
 \end{enumerate}

From this discussion, it follows that it is possible to interpret the degeneracy locus $X$ of 
a general morphism $\varphi :\mathcal{O}^{\oplus m}_{\dP^{n}}\longrightarrow \Omega_{\dP^{n}}(2)$ 
as the set of centres of complexes belonging to a general linear system $\Delta$ of dimension $m-1$ of 
complexes in $\dP^{n}$. Such a $\Delta$ is generated by $m$ independent complexes and corresponds 
to a linear subspace $\dP^{m-1}$ in $\dP(\wedge^{2}V^{\ast})$. The special complexes 
in $\Delta$ are parametrized by the intersection $\dP^{m-1}\cap \check{\mathbb{G}}(1,n)$.
\medskip

We conclude this section by giving some examples.
\begin{itemize}[$\bullet$]
\item  $n=3, m=2$: $X$ is the union of two skew lines;

\item  $n=3, m=3$: $X$ is a smooth quadric, it can be seen as a scroll of degree two over a conic; 

\item  $n=4, m=2$: $X$ is an irreducible conic; 

\item  $n=4, m=3$: $X$ is a smooth projected Veronese surface;

\item  $n=4, m=4$: $X$ is a so-called Segre cubic $3$-fold, which has been exten\-sively studied. It is
singular with singular locus formed by ten distinct points;

\item  $n=5, m=2$: $X$ is the union of three lines in $\dP^{5}$; we shall study the
possible configurations of $X$ in \S 3;

\item  $n=5, m=3$: $X$ is a scroll over a plane cubic, of degree six in $\dP^{5}$, i.e. 
 an elliptic normal scroll in $\dP^{5}$; we shall study this situation in \S 4;

\item  $n=5, m=4$: $X$ is a $3$-fold of degree seven, which is scroll over a cubic
 surface in $\dP^{3}$. It is also known as Palatini scroll (\cite{Ot}).

\end{itemize}

\section{Pencils of linear complexes in $\dP^{5}$}
 
In this section we will study in detail the pencils of linear complexes in $\dP^{5}$, or, 
in other words, maps $\mathcal{O}^{\oplus 2}_{\dP^{5}}\longrightarrow \Omega_{\dP^{5}}(2)$ 
and their degeneracy loci.

We will denote by the same symbol a line both as a subset of $\dP^{5}$ and as a point 
of $\mathbb{G}(1,5)$. Let us recall that the Grassmannian $\mathbb{G}(1,5)$ is an $8$-dimensional variety 
in $\dP^{14}$. A general linear complex of lines in $\dP^{5}$ does not have any centre, 
whereas the special complexes can be of \emph{first or second type}, if they have a line or a $\dP^{3}$ 
as singular space, respectively. The dual space $\check{\dP}^{14}$ parametrizes linear line 
complexes: special complexes correspond to points of the cubic hypersurface $\check{\mathbb{G}}(1,5)
\subseteq \check{\dP}^{14}$. Furthermore special complexes of second type can be interpreted as 
points of the Grassmannian $\mathbb{G}(3,5)$ (which is also embedded in $\check{\dP}^{14}$), 
because a special complex of second type is determined uniquely by its singular space $\dP^{3}$: 
it is formed by the lines intersecting that $\dP^{3}$.

We introduce a rational surjective map:
\[ \psi :\check{\mathbb{G}}(1,5)\dashrightarrow \mathbb{G}(1,5) \]

\noindent which maps a special complex (of first type) to its singular line. So $\psi$ is regular on 
$\check{\mathbb{G}}(1,5)\diagdown \mathbb{G}(3,5)$.
 
The closure of the fibre $\psi^{-1}(\emph{l})$ over a line
\emph{l}
is formed by all special complexes having \emph{l} as singular line. These fibres are linear spaces 
of dimension five; in fact we may interpret $\psi ^{-1}(\emph{l})$ as the linear system of hyperplanes in 
$\dP^{14}$ containing the tangent space to $\mathbb{G}(1,5)$ at the point \emph{l}: 
$T_{\emph{l},\mathbb{G}}\simeq \dP^{8}$. We will use the notation $\dP^{5}_{\emph{l}}$ 
for $\psi ^{-1}(\emph{l})$.

\begin{osservazione} \label{oss1} 
Let $l, m \in \mathbb{G}(1,5)$ be lines of $\dP^{5}$. Then the intersection 
of $\dP^{5}_{l}$ with $\mathbb{G}(3,5)$ is a smooth quadric of dimension \emph{4}. In fact, 
it is formed by the $\dP^{3}$'s containing $l$: this is a Schubert subvariety of the Grassmannian 
and precisely a quadric.

Also the intersection of the fibres $\dP^{5}_{l}\cap \dP^{5}_{m}$ is contained in 
$\mathbb{G}(3,5)$ and is just one point, if $l$ and $m$ do not intersect, or a plane, if they in\-ter\-sect. 
Indeed note that $\dP^{5}_{l}\cap \dP^{5}_{m}$ is the set of the special 
complexes of second type, whose singular space contains both $l$ and $m$. If $l\cap m=\emptyset$, then the 
linear span $\langle l,m \rangle$ is a $\dP^{3}$, so $\dP^{5}_{l}\cap \dP^{5}_{m}$ 
is the point corresponding to the special complex of second type whose singular space is $\langle l,m \rangle$. 
Whereas, if $l\cap m\neq \emptyset$, then $\langle l,m \rangle$ is a plane $\pi$ and $\dP^{5}_{l}
\cap \dP^{5}_{m}$ is the plane corresponding to the linear system of $\dP^{3}$'s containing $\pi$.
\end{osservazione}
                
 Let $L$ be a line in $\check{\dP}^{14}$: $L$ represents a pencil of linear complexes. 
In the general case, $L$ meets $\check{\mathbb{G}}(1,5)$ at three points, which do not belong to 
$\mathbb{G}(3,5)$. Therefore the most general pencil of linear complexes in $\dP^{5}$ contains three 
special complexes of first type. Let $l_{1},l_{2},l_{3}$ be the singular lines of these three complexes. 
We distinguish four possibilities for the reciprocal positions of the three lines.
\medskip 

\textbf{Case 1.} The lines $l_{i}$ are two by two skew and generate the whole $\dP^{5}$. 
This is the most general situation. In this case the three lines do not have any trisecant line.
\smallskip

\textbf{Case 2.} The lines $l_{i}$ are two by two skew and generate a $\dP^{4}$. 
Equivalently, the three lines have one and only one trisecant line.
\smallskip

\textbf{Case 3.} The lines $l_{i}$ are two by two skew and generate a $\dP^{3}$.
\smallskip

\textbf{Case 4.} The lines $l_{i}$ are not two by two skew.
\medskip

Let us analyze separately the four cases. We will see that only in the first case the pencil 
of linear complexes is general.
\smallskip

\textbf{Case 1.} Let $\dP^{5}_{i}:=\psi^{-1}(l_{i})$ be the space that parametrizes special 
complexes having $l_{i}$ as singular line and let $H_{ij}$ be the linear span in $\dP^{5}$ 
of $l_{i}$ and $l_{j}$ for $1\leqslant i<j\leqslant 3$. Note that, by the assumption on the reciprocal 
position of the lines, $H_{ij}\cap H_{ik}=l_{i}$. We will use the same notation for $H_{ij}$ as a $3$-space 
in $\dP^{5}$ and for $H_{ij}$ as a point of $\mathbb{G}(3,5)$. So, by Remark \ref{oss1}, we may write 
$\dP^{5}_{i}\cap\dP^{5}_{j}=H_{ij}$.

We can prove the following: 

\begin{prop} \label{trerette} In the situation of case 1., let $L\subseteq \check{\dP}^{14}$ be a 
line representing a pencil of linear complexes having $l_{1},l_{2},l_{3}$ as singular lines. 
Then $L$ belongs to one of the following families:

\indent a) lines of the plane $\sigma :=\langle H_{ij},1\leqslant i<j\leqslant 3\rangle\subseteq
\check{\dP}^{14}$;

\indent b) lines through $H_{ij}$ intersecting $\dP^{5}_{k}$, for some $1\leqslant i<j\leqslant 3$, 
$k\neq {i,j}$.

In particular, lines $L$ representing general pencils (i. e. pencils not con\-tai\-ning any special complex 
of the second type) correspond to lines of $\sigma$ not passing through any $H_{ij}$.
\end{prop}

\begin{figure}[h]
  \centerline{\epsfig{file=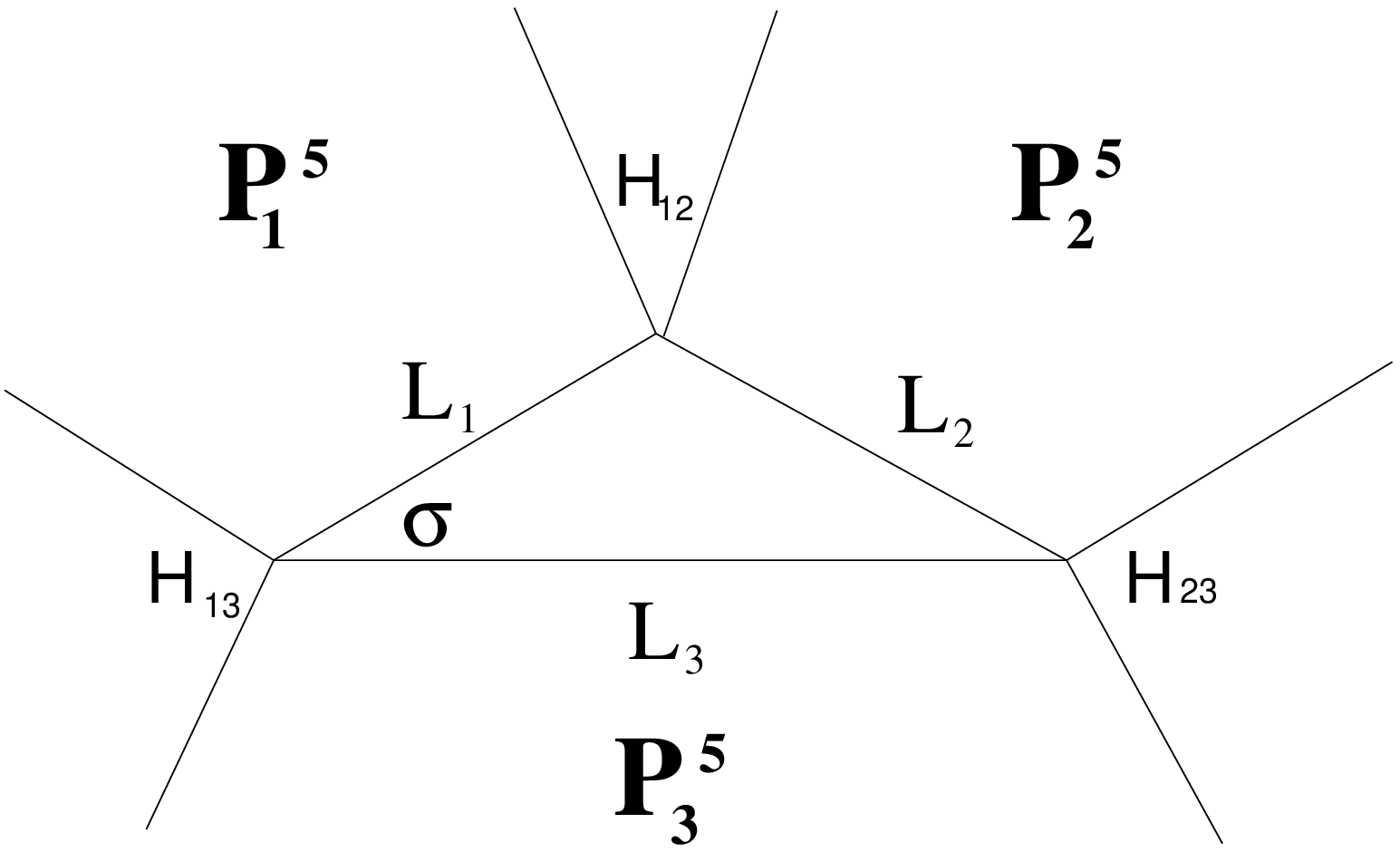,width=26pc}}
  \caption{}
\end{figure}

\emph{Proof}.  We observe first that our assumption on $L$ is equivalent to the fact that $L$ intersects all 
the three spaces $\dP^{5}_{i}$. Moreover $\dP^{5}_{i}\cap \sigma$ is the line $L_{i}$ through 
$H_{ij}$ and $H_{ik}$ (see Figure 1.). So every line in $\sigma$ meets $L_{1},L_{2}$ and $L_{3}$ and 
therefore represents a pencil having $l_{1},l_{2},l_{3}$ as singular lines.

\indent Now let $L\nsubseteq \sigma$ be any line intersecting the three spaces $\dP^{5}_{i}$: 
we want to show that $L$ is a line of type b). We define $R_{i}:=\,L\cap \dP^{5}_{i}$, for $i=1,2,3$. 
Since $L\nsubseteq \sigma$, we can suppose that $R_{2},R_{3}\notin \sigma$. If $R_{2}=R_{3}$, then this 
point is $H_{23}$ and $L$ is of type b). So we assume $R_{2}\neq R_{3}$: then $L$ is contained in 
$\langle \dP^{5}_{2}, \dP^{5}_{3} \rangle$, which is a $\dP^{10}$ by Grassmann relation, 
and $R_{1}=\,L\cap \dP^{5}_{1}$ belongs to $\langle \dP^{5}_{2}, \dP^{5}_{3} \rangle\,
\cap\dP^{5}_{1}=L_{1}$. Now, note that $L$ is contained both in $\langle \dP^{5}_{2}, L_{1}\rangle$ 
and in $\langle \dP^{5}_{3}, L_{1}\rangle$ and therefore in their intersection, which is the plane 
$\sigma$ (again by Grassmann relation). So we have a contradiction that makes us conclude that 
$R_{2}=R_{3}=H_{12}$.

\indent From Remark \ref{oss1}, it follows that each line $L_{i}$ meets $\mathbb{G}(3,5)$ only at the two points 
$H_{ij}$ and $H_{ik}$, otherwise it would be contained in $\mathbb{G}(3,5)$, which is excluded by our 
assumptions. So lines of $\sigma$ not passing through any $H_{ij}$ represent general pencils. 
\ $\Box$

\medskip

\textbf{Case 2.} Assume that $l_{1}, l_{2}, l_{3}$ are two by two skew but generate a hyperplane $H$ in 
$\dP^{5}$. Then they have exactly one trisecant line $r$: it can be constructed as follows. Let $P$ 
be the intersection $\langle l_{1},l_{2}\rangle \cap l_{3}$: $r$ is the only line through $P$ meeting both 
$l_{1}$ and $l_{2}$.

In this case, we can consider the spaces $H_{ij}$, the lines $L_{i}$ and the plane $\sigma$ as in case 
1., but $\sigma$ is now contained in $\mathbb{G}(3,5)$. Indeed, the $H_{ij}$'s intersect two by two along a 
plane in $H$, precisely $H_{ij}\cap H_{ik}=\langle l_{i}, r \rangle$, for all $i,j,k$. This means that the 
lines $L_{1}, L_{2}, L_{3}$, joining the corresponding points in the Pl\"ucker embedding, are completely
contained in $\mathbb{G}(3,5)$, which is a smooth quadric, and this implies that also $\sigma \subseteq 
\mathbb{G}(3,5)$.

The pencils of line complexes having $l_{1}, l_{2}, l_{3}$ as singular lines can be obtained exactly 
as in case 1., but \emph{all of them contain at least one special complex of second type}. In particular, 
those corresponding to the lines in $\sigma$ are pencils of complexes all special of second type.
\medskip

\textbf{Case 3.} If $l_{1}, l_{2}, l_{3}$ are two by two skew but generate a $\dP^{3}$, then 
$H_{12}=H_{13}=H_{23}$ and $\sigma$ ``collapses'' to a point. It is easy to see that a line meeting 
$\dP^{5}_{i}$, for all $i$, necessarily passes through this point, so the corresponding pencil 
contains at least one special complex of second type.
\medskip

\textbf{Case 4.} When at least two among $l_{1}, l_{2}, l_{3}$ intersect each other, several 
configurations are possible. It is clear from the previous discussion that, in all cases, we get pencils 
containing special complexes of the second type.
\medskip

 We can conclude with the following:        

\begin{theorem} If the degeneracy locus $X$ of a map $\varphi : \mathcal{O}^{\oplus 2}_{\dP^{5}}
\longrightarrow \Omega _{\dP^{5}}(2)$ has the expected dimension one, then $X$ is the union of three skew 
lines generating $\dP^{5}$.
\end{theorem}

\emph{Proof}. Indeed, if the pencil of linear complexes corresponding to $\varphi$ contains a special 
complex $\Gamma$ of the second type, then the singular $\dP^{3}$ of $\Gamma$ is contained in $X$. 
So the theorem follows from the previous analysis of cases 1. -- 4.
\ $\Box$ 
\medskip

\begin{osservazione} Let $\mathbb{G}(1,5)^{(3)}$ denote the third symmetric power of the Grassmannian of lines
in $\dP^5$,  which is a projective variety of dimension $24$. $\mathbb{G}(1,\check{\dP}^{14})$, 
of dimension $26$, parametrizes the lines in $\check{\dP}^{14}$, cor\-res\-pon\-ding bi\-jectively to 
pencils of linear line complexes in $\dP^{5}$. There is a natural rational map:
\[ \alpha : \mathbb{G}(1,\check{\dP}^{14})\dashrightarrow \mathbb{G}(1,5)^{(3)}. \]
\noindent $\alpha$ maps a general pencil of linear line complexes of $\dP^{5}$ into the triple of 
its singular lines. The results of this section show that $\alpha$ is dominant and Proposition \ref{trerette} 
describes its general fibre.
\end{osservazione}

\section{Nets of linear complexes in $\dP^{5}$} 

We will study now the nets of linear complexes in $\dP^{5}$, i.e. the maps 
$\mathcal{O}^{\oplus 3}_{\dP^{5}}\longrightarrow \Omega_{\dP^{5}}(2)$ 
and their degeneracy loci.

 As we have seen in \S 2, the singular surface $X$ of a general net $\Delta$ of linear complexes in 
$\dP^{5}$ is an elliptic normal scroll surface, of degree six. Since $\Delta$ is general, it does 
not contain any special complex of second type, so there is a natural surjective regular map, whose fibres 
are lines:
$ \varphi : X\longrightarrow C$,   
where $C$ is the smooth plane cubic, defined by the vanishing of the Pfaffian of the skew-symmetric 
matrix associated to the net. The fibres of $\varphi$ are just the singular lines of the special complexes 
of $\Delta$.    

 We can interpret $\Delta$ as a general plane in $\check{\dP}^{14}$: it intersects 
$\check{\mathbb{G}}(1,5)$ along a cubic curve which is disjoint from $\mathbb{G}(3,5)$. 
Its points represent 
the special complexes of $\Delta$, so the curve can be identified with $C$.

 Let us briefly recall the well-known classification of elliptic normal scrolls in $\dP^{5}$ 
(see \cite{hart}). Every such scroll $X$ is isomorphic to $\dP(\mathcal{E})$, where $\mathcal{E}$ is a rank 
two normalized locally free sheaf over the base curve $C$. The invariant $e$ is necessarily zero in this case, 
so $\deg \mathcal{E}=\,0$. $\dP(\mathcal{E})$ is embedded in $\dP^{5}$ by the very ample linear 
system $|\,C_{0}\,+\,3f\,|$, where $C_{0}$ is a minimal unisecant curve with $C_{0}^{2}=\,0$ and $f$ is a fibre.

\indent There are three cases, according to the Atiyah classification of vector bundles over elliptic 
curves (\cite{atiyah}):
 
\begin{enumerate}
\item $\mathcal{E}=\mathcal{O}_{C}\oplus \mathcal{L}$, with $\mathcal{L}$ a non-trivial
invertible sheaf  of degree $0$. There is a one-dimensional family of such sheaves up to isomorphism, parametrized
by $C$. In this  case, $X$ has two unisecant cubics, $\gamma$ and $\gamma'$. This is the most general case;

\item $\mathcal{E}$ is indecomposable and is unique up to isomorphism. $X$ has only one unisecant cubic;

\item $\mathcal{E}=\mathcal{O}_{C}\oplus \mathcal{O}_{C}$.  $X\simeq \, C\times \dP^{1}$ has a 
$1$-dimensional family of unisecant cubics.
\end{enumerate}      

\par The Hilbert scheme $\mathcal{H}$ of elliptic normal scrolls in $\dP^{5}$ has dimension $36$ and 
is smooth at points representing smooth surfaces (see \cite{Ionescu}). From the above classification, it follows
that  the surfaces of the first type form a dense subset in $\mathcal{H}$, whereas those of the second type form 
a subvariety of dimension $35$ and those of the third type a subvariety of dimension $33$.
  
The nets of linear complexes in $\dP^{5}$ are parametrized by $\mathbb{G}(2,\check{\dP}
^{14})$, the Grassmannian of planes in $\check{\dP}^{14}$, therefore we can define a natural rational 
map:
\[ \rho : \mathbb{G}(2,\check{\dP}^{14})\dashrightarrow \mathcal{H}. \]

\noindent $\rho$ maps a general net $\Delta$ to its singular surface $X$. Since domain and codomain of 
$\rho$ have the same dimension $3$6, the first natural question is if $\rho$ is dominant or, in other words, 
if an elliptic scroll $X$ of the first type is always the degeneracy locus of a suitable bundle map 
$\mathcal{O}^{\oplus 3}_{\dP^{5}}\longrightarrow \Omega_{\dP^{5}}(2)$. In the affirmative case,
 a general fibre of $\rho$ is finite. So the second natural question is what is the degree of that fibre.

\indent An answer to this questions is the main result of this paper. It has been inspired by an article of 
Gino Fano  \cite{Fano}.

\begin{theorem} \label{fano}
Let $X$ be an elliptic normal scroll in $\dP^{5}$ whose associated sheaf is
$\mathcal{E}=
\mathcal{O}_{C}\oplus \mathcal{L}$, with $\mathcal{L}$ non-trivial. Then there are exactly four general nets 
of linear line complexes, whose singular surface is $X$.
\end{theorem}

\indent The proof of Theorem \ref{fano}  will follow after several Lemmas and Re\-marks. We will study first the
relations  between the pro\-per\-ties of a net $\Delta$ having $X$ as sin\-gu\-lar surface and the geometry of
$X$. In particular  we will find some interesting linear series, associated to $\Delta$, on the el\-liptic curve
$C_{X}$ in the  Grass\-man\-nian $\mathbb{G}(1,5)$ re\-pre\-senting the lines of $X$. This will allow us to
explicitly construct  all nets having a fixed $X$ as singular surface. 

\smallskip

 Let $X$ be an elliptic normal scroll of type 1. with unisecant cubics $\gamma$ and $\gamma ^{'}$, 
generating the planes $\pi$ and $\pi ^{'}$ respectively.

\begin{lemma} \label{42}
If $\Delta$ is a general net of linear line complexes in $\dP^{5}$ having $X$ as 
singular surface, then every line contained in $\pi$ or in $\pi ^{'}$ belongs to all complex of $\Delta$.
\end{lemma}

\emph{Proof}. Let $r$ be a general line in $\pi$. Put $r\cap \gamma =\{P_{1},P_{2},P_{3}\}$. Let $l_{i}$ be 
the line of $X$ through $P_{i}$: $l_{i}$ is the singular line of a complex $\Gamma_{i}$ of $\Delta$ and $r$ 
belongs to $\Gamma_{i}$ for $i=$1,2,3. The complexes $\Gamma_{1},\Gamma_{2},\Gamma_{3}$ do not belong to the 
same pencil: indeed the lines $l_{i}$ have the trisecant $r$, so a pencil $\Phi$ containing $\Gamma_{1},
\Gamma_{2},\Gamma_{3}$ should contain at least one special complex of second type (by Theorem \ref{trerette}), but
this  contradicts the assumption that $\Delta$ is general. Therefore $r$ belongs to three independent complexes
of 
$\Delta$, hence to the base locus of $\Delta$. The same conclusion holds true if $r\subseteq \pi ^{'}$.
\ $\Box$  

\begin{lemma} \label{43}
Let $k$ be a line of $X$. The special complexes having $k$ as singular line and containing
all lines  of $\pi$ and $\pi ^{'}$ form a $3$-dimensional linear system. None of these complexes contains all the
lines  of $X$.
\end{lemma}

\emph{Proof}. Let $\dP^{5}_{k}\subseteq \check{\dP}^{14}$ be the space which 
parametrizes the special complexes having $k$ as singular line (notation as in \S 3). The condition that 
the ruled planes $\pi$ and $\pi ^{'}$ are contained in a complex of $\dP^{5}_{k}$ is expressed by 
 two linear conditions, so that these complexes form a $\dP^{3}_{k}\subseteq\dP^{5}_{k}$. 
Indeed, every point of $k$ is a centre of every complex $\Gamma\in\dP^{5}_{k}$, so it suffices to 
impose that one line of $\pi$ (resp. $\pi ^{'}$), not passing through $k\cap \pi $ (resp. $k\cap \pi ^{'}$), 
belongs to $\Gamma$.

 We assume now by contradiction that a complex $\Gamma$ of $\dP^{3}_{k}$ contains all lines 
of $X$. Let us consider the projection $p_{k}$ from $k$ to a complementar space $\dP^{3}$. 
Let $X'\, :=\, p_{k}(X)$: it is a ruled surface of degree $4$. Observe that $r\,:=\,p_{k}(\pi)$ and 
$r'\, :=\,p_{k}(\pi ^{'})$ are skew lines. Every line in $X$, different from $k$, is projected to a line meeting 
$r$ and $r^{'}$. Every point $P\in r$ (resp. $r^{'}$) comes from two points of $\gamma$ (resp. $\gamma ^{'}$), 
which are collinear with $k\cap \gamma $ (resp. $k\cap \gamma ^{'}$). 
Therefore in $X'$ there are two lines through $P$  intersecting $r^{'}$ (resp. $r$).
  
 The projection of $\Gamma$ in $\dP^{3}$ via $p_{k}$ is a linear complex $\Gamma'$ of lines 
of $\dP^{3}$ containing $r, r'$ and all lines of $X'$. If $\Gamma'$ is a singular complex, then it is 
formed by the lines meeting a fixed line, which is impossible because $X'$ is elliptic. If $\Gamma'$ is 
non-singular, the lines of $\Gamma'$ passing through a point $P$ of $r$ should belong to a pencil, 
containing $r$ and the two lines of $X'$ intersecting $r'$: also this is impossible because $r$ and $r'$ 
are skew.    
\ \ $\Box$       

\medskip
Let now $X$ be the singular surface of a net $\Delta$. We will consider the linear series of divisors 
on the elliptic curve $C_{X}$, of degree six in $\mathbb{G}(1,5)$, associated to $X$. The hyperplane 
linear series on $C_{X}$ is a $g^{5}_{6}$: $|\,H_{X}\,|$. If $k$ is a line on $X$, we will   denote 
by $k$ also the corresponding point on $C_{X}$ and  by $\Gamma_{k}$ the line complex of $\Delta$ having 
$k$ as singular line. We can interpret it as a hyperplane in $\dP^{14}$, tangent
 to $\mathbb{G}(1,5)$ at the point corresponding to $k$. Since $\Gamma_{k}$ is also tangent to $C_{X}$ 
at $k$, the intersection of $C_{X}$ and $\Gamma_{k}$, residually to $k$, is a divisor $D_{k}$ of degree four 
on $C_{X}$, such that $D_{k}\,+\,2k$ belongs to the hyperplane series $|\,H_{X}\,|$. If $k$ varies in $C_{X}$,
 the divisors $D_{k}$ describe a non-linear series of degree $4$ and dimension $1$, a $\gamma^{1}_{4}$. 
On the other hand, letting $\Gamma$ vary in $\dP^{3}_{k}$, we obtain a complete linear series 
$|\,D_{k}\,|$, depending on $k$, of dimension three and degree four, i.e. a $g^{3}_{4}$. The following lemma 
is classical.
\medskip

\begin{lemma} \label{EC}
Let $g^{2n-1}_{2n}$ be a complete linear series on an elliptic curve $E$. 
Then there exist exactly four distinct linear series $g^{n-1}_{n}=|\,G_{i}\,|$, i=\emph{1,\ldots,4}, 
such that $|\,2G_{i}\,|=g^{2n-1}_{2n}$.
\end{lemma}

In particular, for all line $k$ there exist four divisors in $|\,D_{k}\,|$ of the form $2E^{(k)}_{i}$,
where $|\,E_{i}\,|$ is a $g^{1}_{2}$.

Similarly there are exactly four linear series such that $|\,H_{X}\,|=|\,2F_{i}\,|$, and $|\,F_{i}\,|$
 is a $g^{2}_{3}$. But $k+E^{(k)}_{i}$ has degree three and $2(k+E^{(k)}_{i})\in |\,2k+D_{k}\,|=|\,H_{X}\,|$. 
Therefore $|\,k+E^{(k)}_{i}\,|$ is one of the four $|\,F_{i}\,|$.

We interpret now $\Delta$ as a plane in $\check{\dP}^{14}$ and its special complexes as the 
points of $\Delta \cap \check{\mathbb{G}}(1,5)$, which is a plane cubic $C$. $C$ and $C_{X}$ are isomorphic 
via the map $\psi$ (see \S 3) which associates to each special complex its singular line.

The hyperplane series $|\,H_{C}\,|$ on $C$ is a $g^{2}_{3}$, representing on $X$ 
triples of lines which are singular loci of pencils contained in $\Delta$.
\medskip
     
\textbf{Claim.} $\psi^{\ast}(H_{X})\in |\,2H_{C}\,|$. Then, up to the isomorphism $\psi$, $|\,H_{C}\,|$ 
is one of the four series $|\,F_{i}\,|$.

\emph{Proof of the claim.} We consider $\psi^{-1}(C_{X} \cap H)$, where $H \in \Delta$ is a hyperplane in 
$\dP^{14}$. Note that a line $k_{\Gamma}$, singular for a complex $\Gamma$, belongs to $H$ if and 
only if $H$ belongs to $\check{k}_{\Gamma}$ in $\check{\dP}^{14}$ and $\check{k}_{\Gamma}$ is 
tangent to $\check{\mathbb{G}}(1,5)$ at all points representing special complexes having $k_{\Gamma}$ as 
singular line, in particular at $\Gamma$. So $\psi^{-1}(C_{X} \cap H)$ is formed by special complexes 
$\Gamma$ of $\Delta$ such that $k_{\Gamma}\in H$, or equivalently that $H\in T_{\Gamma,\check{\mathbb{G}}}$. 
This again is equivalent to the condition that the line $\overline{\Gamma H}$  be tangent to $C=\Delta \cap
\check{\mathbb{G}}$
 at $\Gamma$. So the points of $\psi^{-1}(C_{X} \cap H)$ are in bijection with the tangent lines to $C$ 
passing through $H$. We can conclude that the linear system on $C$ defining $\psi$ contains the system 
cut on $C$ by the polar conics of the points of its plane $\Delta$, which proves the claim.
\medskip

 We have now a rather complete picture of linear series on $C_{X}$, if $\Delta$ is given a priori. In 
particular a $g^{5}_{6}$ and a $g^{2}_{3}$ with $ g^{5}_{6}=2g^{2}_{3}$ are fixed.

We start now from $X$, so the curve $C_{X}$ with hyperplane divisor $H_{X}$ and the four divisors 
$F_{i}$ are fixed. We choose one of them: $F$. Every $k \in C_{X}$ has a residual $g^{1}_{2}$ with respect 
to $F$. By Lemma \ref{EC}, this $g^{1}_{2}$ has four double points $P_{1},\ldots,P_{4}$. 

\medskip
 \textbf{Claim.} There exists a well determined linear complex $\Gamma_{k}$ in 
$\dP^{3}_{k}$ containing the four lines in $X$, which correspond to $P_{1},\ldots,P_{4}$.

\emph{Proof of the claim.} Indeed, if we embed $C_{X}$ in the plane as a plane cubic using $|F|$, then  
$\Gamma_{k}$ corresponds to the conic section cut out by the polar conic of $k$ with respect to $C_{X}$: it
 is tangent to $C_{X}$ at $k$ and it gets through the four contact points of the 
tangent lines to $C_{X}$ passing through $k$. So we have the Claim.                            
 
\medskip
Letting $k$ vary on $C_{X}$, we obtain an elliptic system $S$ of linear complexes $\{\Gamma_{k}\}$. 
To prove Theorem \ref{fano}, it remains to verify that the minimum linear system of complexes containing $S$ is 
a net. Since in $\check{\dP}^{14}$ $S$ is an elliptic curve, it suffices to prove that $S$ is a cubic: 
the plane of $S$ is then the sought net $\Delta$. To compute the degree of $S$, we intersect $S$ with 
$\check{k}$ where $k\in C_{X}$ and $\check{k}$ is the hyperplane in $\check{\dP}^{14}$ which 
parametrizes the complexes containing $k$.

There are only two complexes of $S$ containing $k$, i.e. $\Gamma_{k}$, whose singular line is $k$, and 
$\Gamma ^{'}$, that corresponds to the $g^{1}_{2}$ contained in $|F|$, having $k$ as double element. 
$\Gamma_{k}$ and $\Gamma^{'}$ are the unique points of $\check{k}\cap S$.

 The intersection multiplicity of $\check{k}$ and $S$ at $\Gamma_{k}$ is two, because $\Gamma_{k}$ 
corresponds to a tangent hyperplane to $\mathbb{G}(1,5)$ at $k$. We prove finally that the intersection 
of $\check{k}$ and $S$ at $\Gamma'$ is transversal.
         
 Let $\check{\Gamma'}$ be the hyperplane, through $k$, in $\dP^{14}$, which parametrizes
 the hyperplanes in $\check{\dP}^{14}$ through $\Gamma'$. Observe that 
$\check{\Gamma'}\cap \mathbb{G}(1,5)$ is nothing but the set of the lines of $\Gamma'$. If, by contradiction, 
the intersection of $\check{k}$ and $S$ at $\Gamma'$ is not transversal, then $\check{k}$ is a tangent 
hyperplane to $S$ at $\Gamma'$, so $\check{\Gamma'}$ is tangent to $C_{X}$ at $k$. This is impossible, because 
$k$ is not the singular line of $\Gamma'$. We can conclude that $S$ is just a cubic.

Therefore, we have found a general net of linear complexes in $\dP^{5}$ with $X$ as singular 
surface. There is such a net for each of the four  $g^{2}_{3}$ on $C_{X}$ above described. This completes the 
proof of Theorem \ref{fano}. 
\medskip

We conclude the treatment of this general case with the following Remark, which we will use 
afterwards.  

\begin{osservazione} \label{tipo1}
Let $X$ be an elliptic normal scroll of type 1. with uni\-se\-cant cubics $\gamma$ and $\gamma'$, 
generating the planes $\pi$ and $\pi'$ respectively. The hyperplane series $|H_{\gamma}|$ on $\gamma$ 
and $|H_{\gamma'}|$ on $\gamma'$ are two $g^{2}_{3}$ such that the sum series $|H_{\gamma}+H_{\gamma'}|$ is the 
hyperplane series $|H_{X}|$ on $C_{X}$, up to the Pl\"ucker embedding. 
\end{osservazione}

\emph{Proof}. Observe that a divisor of $|H_{\gamma}+H_{\gamma'}|$ corresponds to the union of $r\cap
\gamma$ and $r'\cap\gamma'$, where $r\subseteq\pi$ and $r'\subseteq\pi'$ are two skew lines. The 
$3$ -space $\langle r,r'\rangle$ determines uniquely the special complex $\Gamma$ of second type, 
formed by all and only lines of $\dP^{5}$ meeting $\langle r,r'\rangle$. $\Gamma$ corresponds 
to a hyperplane section of $C_{X}$, and precisely to a divisor of $|H_{X}|$. This proves the Remark.
\ $\Box$   
\medskip

We may use the same arguments, with suitable changes, to study also elliptic scrolls $X$ of 
type 2. and 3., but in both cases the situation results to be different.

Let $X$ be an elliptic scroll of type 2., i.e. with  unisecant cubic $\gamma$ generating the plane
$\pi$. Such a scroll has the property that the $4$-space, generated by three lines of $X$ meeting $\gamma$ at
collinear  points, always contains the plane $\pi$. This follows immediately from the fact that $X$ is embedded
in 
$\dP^{5}$ by the very ample linear system $|\,C_{0}+3f\,|$, where $C_{0}$ is the unique unisecant 
cubic and $f$ is a fibre, as recalled at the beginning of this section.

\begin{lemma}
If $\Delta$ is a general net of linear line complexes in $\dP^{5}$ having $X$ as singular surface, 
then every line contained in $\pi$ belongs to all complex of $\Delta$.
\end{lemma}

\emph{Proof}.
See the proof of Lemma \ref{42}.
\ $\Box$

\begin{lemma}\label{46}
Let $k$ be a line of $X$. The special complexes having $k$ as singular line and containing all lines of $\pi$ 
form a $4$-dimensional linear system.
Only one of this complexes contains all the lines of $X$.
\end{lemma}

\emph{Proof}. 
For the first assertion see the proof of Lemma \ref{43}.
Let us consider the projection $p_{k}$ from a line $k$ of $X$ to a complementar space $\dP^{3}$.
Let $X':=p_{k}(X)\,$: it is a ruled surface of degree $4$ with only one unisecant line $r:=p_{k}(\pi)$. Observe 
that every line in $X$, different from $k$, is projected to a line meeting $r$, and every point $P\in r$ 
comes from two points of $\gamma$, which are collinear with $k\cap \gamma$. Therefore in $X'$ there are two 
lines through $P$ intersecting $r$.

Observe that the lines of $X'$ are contained in the special complex $\Gamma'$ of lines of 
$\dP^{3}$ with $r$ as singular line. Let $\Gamma :=p_{k}^{-1}(\Gamma')\,$: it is a linear complex of 
$\dP^{5}$, containing all the lines of $X$ and, moreover, having $k$ as singular line. Indeed, $\Gamma'$ 
is formed by the union of the (closure of the) fibres $p_{k}^{-1}(l)$, when $l$ varies among the lines meeting 
$r$. Note that $p_{k}^{-1}(l)$ is the set of the lines contained in the $3$-space $\langle\,k,l\,\rangle$.
\ $\Box$

\begin{osservazione}
If $X$ is an elliptic scroll of type 2., then, for all line $k$ in $X$, there is a unique linear special complex 
$H_{k}$ in $\dP^{5}$, having $k$ as singular line and containing all the lines of $X$. So we obtain a 
$1$-dimensional elliptic system $\{H_{k}\}_{k\subseteq X}$ of special complexes containing $X$.         
\end{osservazione}

Now, let $X$ be of type 3., i.e. $X\simeq\gamma\times\dP^{1}$, where $\gamma$ is a cubic 
generating the plane $\pi$.

\begin{lemma}
If $\Delta$ is a general net of linear line complexes in $\dP^{5}$ having $X$ as singular surface, 
then the base locus of $\Delta$ contains all the $\infty^{1}$ ruled planes of the cubics in $X$.
\end{lemma}
\emph{Proof}.
See the proof of Lemma \ref{42}.
\ $\Box$ 
\medskip

\begin{lemma}
Let $k$ be a line of $X$. Among the special complexes having $k$ as singular line, there is a linear system
of dimension $2$ of such complexes containing all the lines of $X$.
\end{lemma}

\emph{Proof}.
Let us consider the projection $p_{k}$ from a line $k$ of $X$ to a complementar space $\dP^{3}$. 
$X':=p_{k}(X)$ is a quadric surface and $p_{k}$ is a $2:1$ map.

Observe that there is a linear system of dimension $2$ of complexes of lines of $\dP^{3}$ 
containing the lines of the rulings, projection via $p_{k}$ of the lines of $X$ on $X'$. These complexes 
can be lifted to special complexes of $\dP^{5}$ containing $X$ and having $k$ as 
singular line. To prove this, we may use the same arguments of the proof of Lemma \ref{46}.
\ $\Box$ 
\medskip 

If we start from $X$ of type 2. or 3., i.e. the curve $C_{X}$ with hyperplane divisor $H_{X}$ and the 
four series $|F_{i}|$, such that $|2F_{i}|=|H_{X}|$, for $i=1,\ldots,4$, are fixed, we find that \emph{one} of 
the four nets of complexes, having $X$ as singular surface, contains some special complexes of second type. 
Indeed we have the following:

\begin{osservazione}
Let $X$ be an elliptic normal scroll of type 2. or 3. The two series $g^{2}_{3}$ defined in Remark  \ref{tipo1}
coincide, because of the geometry of $X$. If  $|H_{\gamma}|$ denotes this $g^{2}_{3}$, then  
$|2H_{\gamma}|=|H_{X}|$ and so $|H_{\gamma}|$ is one of the four $F_{i}$.
\end{osservazione}

The net $\Delta$ of complexes having $X$ as singular surface, that can be constructed starting from 
the series $|H_{\gamma}|$, contains complexes of second type. Indeed, for every divisor of 
$|H_{\gamma}|$, we obtain three lines $l_{1},l_{2},l_{3}$ of $X$ having a common trisecant, so that among 
the complexes of the pencil contained in $\Delta$ and having $l_{1},l_{2},l_{3}$ as singular lines, there are 
certainly special complexes of second type (see \S 3).

{} From this discussion we conclude that the singular locus of such a net $\Delta$ has codimension $2$
instead of the  expected codimension $3$. Indeed, it is  the union of $X$ with the singular $3$-space
of the  special complexes of second type of $\Delta$.

\end{article}
\end{document}